\documentclass[a4paper,11pt]{article}
\usepackage[ansinew]{inputenc}
\usepackage[english]{babel}
\usepackage{mathrsfs}
\usepackage{amsfonts}
\usepackage[centertags]{amsmath}
\usepackage{amsthm}
\usepackage{latexsym,amsmath,amssymb}
\usepackage{graphicx,color}
\usepackage{eqnarray}
\usepackage[all]{xy}
\input diagxy
\newtheorem{teo}{Theorem}[section]
\newtheorem{prop}[teo]{Proposition}
\newtheorem{defin}[teo]{Definition}

\newtheorem{state}[teo]{Problem}

\newtheoremstyle{obs}
  {3pt}
  {3pt}
  {}
  {}
  {\bfseries}
  {.}
  {.5em}
  {}
\theoremstyle{obs}

\newtheorem{remark}[teo]{Remark}

\newcommand{\ds}{\displaystyle}

\def\parx{\frac{\partial}{\partial x}}
\def\pary{\frac{\partial}{\partial y}}
\def\parz{\frac{\partial}{\partial z}}
\def\parvx{\frac{\partial}{\partial v_x}}
\def\parvy{\frac{\partial}{\partial v_y}}
\def\parvz{\frac{\partial}{\partial v_z}}

\def\tabaddress#1{{\small\it\begin{tabular}[t]{c}#1
\\[1.2ex]\end{tabular}}}

\title{Constraint Algorithm for Extremals in Optimal Control Problems}
\author{{\sc Mar\'ia Barbero-Li\~n\'an}
\thanks{{\bf e}-{\it mail}: mbarbero@ma4.upc.edu}, {\sc Miguel C. Mu\~noz-Lecanda\thanks{{\bf e}-{\it mail}:
matmcml@ma4.upc.edu}}
\\
 \tabaddress{Departamento de Matem\'atica Aplicada IV\\
  Edificio C-3, Campus Norte UPC.
  C/ Jordi Girona 1. E-08034 Barcelona. Spain}
}

\date{January 29, 2008}

\begin{document}
\maketitle \pagestyle{myheadings}

\thispagestyle{empty}

\begin{abstract}
A geometric method is described to characterize the
different kinds of extremals in optimal control theory.
This comes from the use of a presymplectic constraint
algorithm starting from the necessary conditions given by
Pontryagin's Maximum Principle. Apart from the design of
this general algorithm useful for any optimal control
problem, it is showed how it works to split the set of
extremals and, in particular, to characterize the strict
abnormality. An example of strict abnormal extremal for a
particular control-affine system is also given.
\end{abstract}

  \bigskip
  {\bf Key words}:  Pontryagin's Maximum Principle, extremals, optimal
  control problems, abnormality, strict abnormality,
  presymplectic.
\bigskip

\vbox{\raggedleft AMS s.\,c.\,(2000): 34A26,  49J15, 49K15,
70G45, 70H05, 70H45 }\null

\section{Introduction}\label{intro}

A difficult problem in optimal control is to obtain
extremals, that is, curves candidates to be optimal
solutions. The usual way to deal with that is through
successive differentiations of some necessary conditions
for optimality, see for instance
\cite{2004Agrachev,BookBC,1977Krener}. Here we not only
give a method to split the set of extremals for any optimal
control problem, but also explain geometrically the meaning
of the successive differentiations.

There are different kinds of extremals: normal, abnormal
and strictly abnormal. The abnormal extremals have been
partially ignored for years. In the nineties, the papers by
R. Montgomery, W. Liu and H. J. Sussmann \cite{M94,LS96}
showed up the importance of analyzing abnormal extremals,
because they can be optimal. Therefore, the search for
abnormal and also strict abnormal extremals has become an
appealing issue for the last fourteen years, as is backed
by
\cite{AgrachevGauthier,AZ,BookBC,2006BonnardTrelat,ChitourJeanTrelatDiffGeom2006,
ChitourJeanTrelat2006,Trelat2000,Trelat2006}. The main
attraction of abnormality is its exclusive dependence on
the geometry of the control system and for the strict
abnormality is the fact that the strict abnormal extremals
are not normal.

The essential result to describe the general method of this
paper and to have techniques to solve optimal control
problems is Pontryagin's Maximum Principle, despite only
providing necessary conditions for optimality. Although the
natural geometric framework of Pontryagin's Maximum
Principle is the symplectic one
\cite{1997Jurdjevic,2006AndrewCourse,S98Free}, to our
purpose the presymplectic formalism will be more useful
\cite{MarinaIbortPanorama,Eduardo04}. Then, we have an
implicit equation including some compatibility conditions,
that must be satisfied in order to have solution, besides
the dynamical Hamilton's equations. The former is a
necessary condition of the maximization of the hamiltonian
over the controls according to the classic Maximum
Principle \cite{67LeeMarkus,P62}. Hence, in the
presymplectic framework a weaker version of Pontryagin's
Maximum Principle is stated. Instead of the above classic
necessary condition, we have an implicit differential
equation that sets up a constraint algorithm in the sense
given in \cite{Pepin,1979Gotay,1978Gotay,1992GraciaPons}.
This presymplectic algorithm comes from the Dirac-Bergmann
theory of constraints developed in the fifties for quantum
field theory. This algorithm has been already adapted and
used to study singular optimal control problems
\cite{MarinaIbortPanorama} and to study optimal control
problems with nonholonomic constraints \cite{2004Madrid}.

The aim of this paper, according to the previous optimal
control formulation, is to give a precise and geometric
description of how to use the constraint algorithm to
determine where the dynamics of normal extremals takes
place and also the dynamics of abnormal ones. We also
obtain sufficient conditions to have both kinds of
extremals. These conditions elucidate how to determine the
strict abnormality. This adaptation of the algorithm to the
study of the extremals is mostly developed in \S
\ref{Scharacterize}, under the assumption of the control
set being open and the differentiability with respect to
the controls whenever is needed.

The importance and the generality of the theory elaborated
can be highlighted by the revisit of the characterization
of abnormal extremals in some known examples such as
subRiemannian geometry and single-input control-affine
systems. Using the algorithm and distinguishing different
cases that come up, it may be checked that some of the
situations correspond with the results obtained by A. A.
Agrachev, Y. Sachkov, I. Zelenko, W. Liu and H.J. Sussmann
\cite{2004Agrachev,AZ,LS96}. Our method collects all their
results on the existence of abnormal extremals. So the
described procedure allows us to study geometric
 and generically the extremals for any control system and
obtain the dynamics of the abnormal extremals in a natural
and understandable way.

The organization of the paper is as follows: in Section 2
after a brief review of some notions in optimal control
theory, we state the optimal control problem and
Pontryagin's Maximum Principle in the suitable framework
for this paper, that is, in the presymplectic one. Section
3 concentrates on the new material, so it is devoted to
describe the geometric process used to characterize
extremals in optimal control problems with fixed time.
After studying the fixed time problem, we explain how the
algorithm works for the free time case in Section 4.
Finally, in Section 5, we find a strict abnormal extremal
for a control-affine system using the presymplectic
constraint algorithm.

In the sequel, all the manifolds are real, second countable
and ${\cal C}^{\infty}$. The maps are assumed to be ${\cal
C}^{\infty}$. Sum over repeated indices is understood.

\section{Presymplectic optimal control problems}\label{SOCP}

A {\sl control system} is defined by a set of differential equations
depending on parameters. More precisely, let $M$ be a smooth
manifold, ${\rm dim} \, M=m$, $U$ be an open set of $\mathbb{R}^k$
called the {\sl control set} with $k\leq m$. A vector field $X$
along the projection $\pi\colon M\times U \rightarrow M$ is a map
$X\colon M\times U \rightarrow TM$ such that the following diagram
is commutative
$$\bfig\xymatrix{&TM\ar[d]_{\txt{\small{$\tau_M$}}}\\
M\times U \ar[ur]
^{\txt{\small{$X$}}}\ar[r]^{\txt{\small{$\pi$}}} &M}\efig$$
where $\tau_M$ is the natural projection of the tangent
bundle. We denote the set of these vector fields as
$\mathfrak{X}(\pi)$. A {\sl control system} is an element
of $\mathfrak{X}(\pi)$.

Let $I\subset \mathbb{R}$, a curve $(\gamma,u)\colon I
\rightarrow M\times U$ is an integral curve of $X$ if
\begin{equation}\label{eqX} \dot{\gamma}=X\circ
(\gamma,u), \quad {\rm that} \quad {\rm is},  \;\;
\dot{\gamma}(t)=X(\gamma(t),u(t)).\end{equation} Now, we
can introduce a cost function $F \colon M\times U
\rightarrow \mathbb{R}$ and the functional
$${\cal S}[\gamma,u]=\int_I F(\gamma,u)\, dt$$
defined on curves $(\gamma, u)$ with a compact interval as
domain. We are interested in the following problem:
\begin{state} \textbf{(Optimal Control Problem,
$OCP$)}\label{OCP}
\\
Given the elements $M$, $U$, $X$, $F$, $I=[a,b]$, $x_a$,
$x_b\in M$. Find $(\gamma,u)$ such that
\begin{itemize}
\item[(1)] the endpoint conditions are satisfied $\gamma(a)=x_a$, $\gamma(b)=x_b$,
\item[(2)] $ \dot{\gamma}(t)=X(\gamma(t), u(t))$, $t\in I$, and
\item[(3)] ${\cal S}[\gamma,u]$ is minimum over all curves on $M\times U$
satisfying  $(1)$ and $(2)$.
\end{itemize}
\end{state}
A solution $(\gamma,u)$ to this problem is called {\sl
optimal curve}. The mappings $(\gamma,u) \colon I
\rightarrow M\times U$ are piecewise differentiable and the
vector field $X$ along $\pi$ and the cost function $F
\colon M \times U \rightarrow \mathbb{R}$ are
differentiable enough.

\subsection{Pontryagin's Maximum Principle}\label{SPMP}

As was said in \S \ref{intro}, we state Pontryagin's
Maximum Principle from a presymplectic viewpoint
\cite{MarinaIbortPanorama,2003BcnControlSim,2004Madrid,Eduardo04}.
In this approach, the main elements are:
\begin{itemize}
\item The {\sl presymplectic manifold} $(T^*M\times
U,\Omega)$, where $\Omega$ is the closed 2-form
 on $T^*M\times U$ given by the pull-back through
$\pi_1\colon T^*M\times U\rightarrow T^*M$ of the canonical
2-form on $T^*M$.
\item A {\sl presymplectic Hamiltonian system}
$(T^*M\times U, \Omega, H)$, where $H\colon T^*M\times U
\rightarrow \mathbb{R}$ is the {\sl Pontryagin's
Hamiltonian function} given by
\begin{equation*}
H(\lambda,u)=\langle \lambda, X(x,u)\rangle + p_0
F(x,u)=H_X(\lambda,u) + p_0 F(x,u),
\end{equation*} with $\lambda\in T^*_xM$, $p_0\in
\{-1,0\}$ and the notation $H_X(\lambda,u)=\langle \lambda,
X(x,u)\rangle$.
\end{itemize}
Observe that the kernel of $\Omega$ contains the
$\pi_1$-vertical vector fields, that is,
$\pi_1$-projectable vector fields $Z\in
\mathfrak{X}(T^*M\times U)$ such that $(\pi_1)_*Z=0$. Thus,
$\Omega$ is degenerate. For details in presymplectic
formalism see
\cite{Pepin,1979Gotay,1978Gotay,1992GraciaPons,Eduardo04,92Bcn}.

\begin{teo} \textbf{(Pontryagin's Maximum
Principle, presymplectic form)}\label{PMPpre} Let $U\subset
\mathbb{R}^k $ be an open set and $(\gamma,u)\colon [a,b]
\rightarrow M \times U$ be a solution of the optimal
control problem \ref{OCP} with endpoint conditions $x_a$,
$x_b$. Then there exist $\lambda\colon [a,b] \rightarrow
T^*M$ along $\gamma$ (i.e. the natural projection of
$\lambda$ to $M$ is $\gamma$), and a constant $p_0\in
\{-1,0\}$ such that:
\begin{enumerate}
\item $(\lambda,u)$ is an integral curve of a
Hamiltonian vector field $X_H$ that satisfies
\begin{equation}\label{preqP}
i_{X_H}\Omega={\rm d}H, \; {\rm that \; is}, \;
i_{(\dot{\lambda}(t),\dot{u}(t))}\Omega={\rm d}H(\lambda(t),u(t));
\end{equation}
\item \begin{itemize}
\item[(a)] $H(\lambda(t),u(t))$ is constant everywhere in $t\in [a,b]$;
\item[(b)] $(p_0,\lambda(t))\neq 0$ for each  $t \in
[a,b]$.
\end{itemize}
\end{enumerate}
\end{teo}
As $\Omega$ is degenerate, (\ref{preqP}) does not have
solution in the whole manifold $T^*M\times U$. As explained
in $\S$ \ref{Scharacterize}, it may have a solution if we
restrict the equation to the submanifold defined implicitly
by
\[S=\{\beta\in T^*M\times U \, | \, i_v \, {\rm d}H=0, \quad {\rm for }\; v\in \ker
\Omega_{\beta} \},\] and locally, $S=\{\beta\in T^*M\times
U \, | \, \ds{\frac{\partial H}{\partial u_l}(\beta)=0} ,
\quad l=1,\ldots, k\}$.

\remark \label{aclara} Observe that this is a necessary
condition for the Hamiltonian to have an extremum over the
controls as long as $U$ is an open set. In the classic
Pontryagin's Maximum Principle \cite{P62}, the Hamiltonian
is equal to the maximum of the Hamiltonian over the
controls. Therefore, Theorem \ref{PMPpre} is a weaker
version of the classic Maximum Principle.

The necessary conditions 1-2 of Theorem \ref{PMPpre}
determine different kinds of extremals.

\begin{defin}\label{definextremal}
A curve $(\gamma, u)\colon [a,b] \rightarrow M \times U$ is
\begin{enumerate}
\item an \textbf{extremal for $OCP$} if there exist $\lambda\colon [a,b] \rightarrow
T^*M$ and a constant $p_0\in \{-1,0\}$ such that $(\lambda,
u)$ satisfies the necessary conditions of Pontryagin's
Maximum Principle;
\item a \textbf{normal extremal for $OCP$}
if it is an extremal and $p_0=-1$, that is, the Hamiltonian
is $H^{[-1]}=H_X - F$;
\item an \textbf{abnormal extremal for $OCP$} if it is an extremal and
$p_0=0$, that is, the Hamiltonian is $H^{[0]}=H_X$;
\item a \textbf{strictly abnormal
extremal for $OCP$} if it is not a normal extremal, but it
is an abnormal extremal;
\end{enumerate}
The curve $(\lambda,u)\colon [a,b]\rightarrow T^*M\times U$
is called \textbf{biextremal for $OCP$}.
\end{defin}
Pontryagin's Maximum Principle lifts optimal solutions to the
cotangent bundle. The uniqueness of the lifts is not guaranteed,
that is, some extremals could be lifted in two different ways:
normal and abnormal.

\section{Characterization of extremals}\label{Scharacterize}

Here we take advantage of the necessary conditions in
Theorem \ref{PMPpre} to determine where the different kinds
of extremals above defined are contained. We are specially
interested in strict abnormal extremals and abnormal
extremals as a consequence of \cite{LS96,M94}. A meaningful
and constructive procedure in presymplectic manifolds in
order to find a solution to Problem \ref{presymplProblem}
is the constraint algorithm
\cite{Pepin,1979Gotay,1978Gotay,1992GraciaPons,92Bcn}.
\begin{state}\label{presymplProblem} Given a presymplectic
system $(M,\Omega, H)$, find $(N,X)$ such that
\begin{itemize}
\item[(a)] $N$ is a submanifold of $M$,
\item[(b)] $X \in \mathfrak{X}(M)$ is tangent to $N$ and verifies $i_X\Omega= {\rm d}H$ on $N$,
\item[(c)]$N$ is maximal among all the submanifolds
satisfying $(a)$ and $(b)$.
\end{itemize}
\end{state}
As mentioned in $\S$ \ref{SPMP} the presymplectic equation
(\ref{preqP}), $i_{X_H}\Omega={\rm d}H$, has solution in the primary
constraint submanifold $N_0=\{x\in M \, | \, \exists \, v_x \in T_xM
\, , \quad i_{v_x}\Omega={\rm d}_xH\}$, or equivalently, $N_0=\{x\in
M \, | \, ({\rm L}_Z H)_x=0 \, , \; \forall \; Z\in \ker \Omega\}$,
where ${\rm L}_Z$ is the Lie derivative with respect to $Z$. See
\cite{1978Gotay,92Bcn} for details on these equivalence.

Locally, $ X_H=A^i\partial/\partial x^i
+B_j\partial/\partial p_j+C_l\partial/\partial u_l$ with
$A^i=\partial H/\partial p_i$ and $B_j=-\partial H/\partial
x^j$ because of the presymplectic equation (\ref{preqP}).
Moreover, as $\partial / \partial u_l \in \ker \Omega$,
\begin{equation}\label{N0} N_0=\{(\lambda,u)\in T^*M\times U \, | \,
\frac{\partial H}{\partial u_l}=\lambda_j \, \frac{\partial
X^j}{\partial u_l}+p_0\frac{\partial F}{\partial u_l}=0 \;
, \quad l=1,\ldots,k\}.
\end{equation} The solution on $N_0$ is not
necessarily unique. Indeed, if $X_0$ is a solution, then $X_0+\ker
\Omega$ is the set of all the solutions. We may consider $X_0$ as a
vector field defined on the whole $M$ because $N_0$ is closed and we
assume that $N_0$ is a submanifold of $M$.

Take the pair $(N_0, X_0+\ker \Omega)$, rewritten as
$(N_0,X^{N_0})$. Observe that we are looking for an element in
$X^{N_0}$ tangent to $N_0$. Then,
\[N_1=\{x\in N_0 \, | \, \exists \, X\in X^{N_0} \, ,  \quad  X(x)\in
T_xN_0\},\] locally
\begin{eqnarray}\label{N1} N_1=\{(\lambda,u)\in N_0 & | & 0=X_H\left(\frac{\partial
H}{\partial u_l}\right) =\partial H/\partial p_i\partial^2
H/\partial x^i\partial u_l -\partial H/\partial
x^j\partial^2H/\partial p_j\partial u_l\nonumber \\
&&+C_r\partial^2H/\partial u_r\partial u_l , \quad
l=1,\ldots, k\}.
\end{eqnarray}
 If the matrix $(\partial
H/\partial u_r \partial u_l)_{rl}$, multiplying $C_r$, is
not invertible, the OCP is singular
\cite{MarinaIbortPanorama}, otherwise it is regular.

This step stabilizes the constraints in $N_0$ providing a
new pair $(N_1,X^{N_1})$ where $X^{N_1}$ is the set of the
vector fields solution and tangent to $N_0$. Inductively,
we arrive at $(N_i,X^{N_i})$ where we assume that $N_i$ is
a submanifold of $M$ and we define $N_{i+1}=\{x\in N_i \, |
\, \exists \, X\in X^{N_i} \, , \quad X(x)\in T_xN_i\}$,
obtaining the sequence
\[M \supseteq N_0 \supseteq N_1 \supseteq \ldots \supseteq
N_i \supseteq N_{i+1} \supseteq \ldots\] and the corresponding
$X^{N_{i+1}}$. Let
\[N_f=\bigcap_{i \geq 0}N_i, \quad X^{N_f}=\bigcap_{i \geq
0}X^{N_i},\] if $(N_f,X^{N_f})$ is a nontrivial pair, it is the
solution to Problem \ref{presymplProblem}. If at one step
$N_i=N_{i+1}$, the algorithm finishes with $N_f=N_i$.

Note that each step of the algorithm can reduce the set of
points of $M$ where there exists solution, that is,
$N_i\subsetneq N_{i-1}$, and can also reduce the degrees of
freedom of the set of vector fields solution,
$X^{N_i}\subsetneq X^{N_{i-1}}$. In terms of control
systems, the desirable objectives are to restrict the
problem to a smaller submanifold of $T^*M\times U$ and to
determine the input controls. Observe that, generally, a
step of the algorithm can provide us new constraints and
the determination of some controls at the same time. Hence,
either a unique vector field is found or the new
constraints must be stabilized or the set must be split in
submanifolds. At the final step, we have either a unique or
nonunique vector field and a submanifold that could be an
empty or discrete set.

 \remark Observe that we do not miss any extremal using the
constraint algorithm, in contrast to what happens in subRiemannian
geometry in \cite{LS96}, where using a less geometric approach they
miss the constant extremals.

Now, let us focus again on optimal control problems where
there are two distinct Hamiltonians depending on the value
of the constant $p_0$. Thus, from (\ref{N0}) it is deduced
that the constraint algorithm must be run twice, one for
each Hamiltonian, as is explained in \S \ref{characabn},
\ref{characnorm}.

\subsection{Characterization of abnormality}\label{characabn}

First, we characterize a subset of $T^*M\times U$ where the abnormal
biextremals are, if they exist. In this situation $p_0=0$ and the
corresponding Pontryagin's Hamiltonian is $H^{[0]}=H_X$. Then the
primary constraint submanifold (\ref{N0}) becomes
\begin{equation}\label{N00} N_0^{[0]}=\{(\lambda,u) \in
T^*M \times U \, | \, \lambda_j \, \frac{\partial
X^j}{\partial u_l}=0, \quad l=1,\ldots, k\},\end{equation}
the submanifold (\ref{N1}) is
\begin{equation*} N_1^{[0]}=\{(\lambda,u) \in
N_0^{[0]} \, | \, \lambda_j \, ( X^i\frac{\partial^2
X^j}{\partial x^i\partial u_l}-\frac{\partial X^j}{\partial
x^i}\frac{\partial X^i}{\partial u_l}+C_r\frac{\partial^2
X^j}{\partial u_r\partial u_l})=0, \quad l=1,\ldots,
k\},\end{equation*} and the algorithm continues.

Once we have the final constraint submanifold $N_f^{[0]}$ for
abnormality, we have to delete the biextremals in the zero fiber
because these biextremals do not satisfy the necessary condition
(2.b) of Pontryagin's Maximum Principle \ref{PMPpre}. For the sake
of simplicity and clarity, we denote this actual final constraint
submanifold with the same name $N_f^{[0]}$.

\begin{prop} If $N^{[0]}_f\neq \emptyset$, that is, $(\lambda,u)\in N^{[0]}_f$
with $\lambda\neq 0$, then $(\gamma,u)=(\pi_M\times {\rm
Id})(\lambda,u)$ is an abnormal extremal.
\end{prop}
\subsection{Characterization of normality}\label{characnorm}

Analogous to $\S$ \ref{characabn}, for $p_0=-1$,
Pontryagin's Hamiltonian is $H^{[-1]}=H_X-F$.

Then the primary constraint submanifold (\ref{N0}) becomes
\begin{equation}\label{N01} N_0^{[-1]}=\{(\lambda,u) \in
T^*M \times U \, | \, \lambda_j \, \frac{\partial
X^j}{\partial u_l}-\frac{\partial F}{\partial u_l}=0, \quad
l=1,\ldots,k\},\end{equation} the submanifold (\ref{N1}) is
\begin{equation*} \begin{array}
{ll} N_1^{[-1]}=\{(\lambda,u) \in N_0^{[-1]} \, | \, &
\ds{\lambda_j \, (X^i\frac{\partial^2 X^j}{\partial
x^i\partial u_l}-\frac{\partial X^j}{\partial
x^i}\frac{\partial X^i}{\partial u_l})-X^i\frac{\partial^2
F}{\partial x^i\partial u_l}}  \\  & \\
& \ds{+C_r(\lambda_j\frac{\partial^2 X^j}{\partial
u_r\partial u_l}-\frac{\partial^2 F}{\partial u_r\partial
u_l})=0, \quad l=1,\ldots, k}\},\end{array}\end{equation*}

In general, the determination of the controls for normal
extremals depends on the given cost function, unless it is
quadratic or linear or independent on the controls. To a
better understanding of all this process we address the
reader to the examples in $\S$ \ref{SFOCP},
\ref{SExamples}.

It can be observed that Hamilton's equations for $\dot{x}^i$ are the
same for both Hamiltonian functions, for $p_0=0$ and $p_0=-1$, since
the cost function does not depend on the momenta $p$'s. Hamilton's
equations for $\dot{p}_i$ are equal for cost functions not depending
on $x^i$. For instance, if the cost function is constant, as in the
case of time-optimal.

The final constraint submanifolds $N^{[0]}_f$ and $N^{[-1]}_f$
restrict the set of points where the biextremals of the Optimal
Control Problem \ref{OCP} are. But, even in the case that Hamilton's
equations are the same, $N^{[0]}_f$ and $N^{[-1]}_f$ could be
different. Then the integral curves of the same vector field in
$T^*M\times U$ along the same extremal in $M$ may be different
depending on where the initial conditions for the momenta are taken.
In other words, there may exist abnormal extremals being normal and
viceversa. For a deeper study about how the extremals are we need to
project the biextremals on the base manifold $M\times U$ using
$\rho_1=\pi_M \times {\rm Id}\colon T^*M\times U \rightarrow M\times
U$.

Summarizing all the above comments, we have the following
propositions.
\begin{prop}\label{c1} If there exists $(\lambda,u)\in N^{[-1]}_f$,
then $(\gamma,u)=(\pi_M\times {\rm Id})(\lambda,u)$ is a
normal extremal.
\end{prop}

\begin{prop}\label{c2} Let $(\gamma,u)$ be an abnormal extremal. If there exists a covector $\lambda$ along $\gamma$
such that $(\lambda,u)\in N^{[-1]}_f$, then $(\gamma,u)$ is
also a normal extremal.

Let $(\gamma,u)$ be a normal extremal. If there exists a
covector $\lambda$ along $\gamma$ such that $(\lambda,u)\in
N^{[0]}_f$, then $(\gamma,u)$ is also an abnormal extremal.
\end{prop}

\begin{prop}\label{c3} If there exist $(\lambda^{[0]},u^{[0]})\in N^{[0]}_f$
with $\lambda^{[0]}\neq 0$ and $(\lambda^{[-1]},u^{[-1]})\in
N^{[-1]}_f$ such that
$\pi_M(\lambda^{[0]})=\pi_M(\lambda^{[-1]})=\gamma$, then $\gamma$
is an abnormal extremal and also a normal extremal.
\end{prop}

\remark \label{nocontrol} In this last proposition we do not
consider the control as a part of the extremal, because it may
happen that different controls give the same extremals in $M$
depending on the control system. So we project onto $M$ the
biextremals to compare them. Under some assumptions about the
control systems, such as control-affinity with independent control
vector fields, different controls give different extremals. If so
happens, we will project the biextremals onto $M\times U$ through
$\rho_1$ to compare them.

\remark \label{setextrem} The union of both final
constraint submanifods do not cover exactly the set of
extremals in Definition \ref{definextremal}, because the
condition $(2.a)$ in Theorem \ref{PMPpre} is not included
in the final constraint submanifolds. See $\S$ \ref{SFOCP}
to get a better understanding.

\subsection{Characterization of strict
abnormality}\label{characstrict}

The normal and abnormal extremals in Definition \ref{definextremal}
do not constitute a disjoint partition of the set of extremals as
propositions \ref{c2}, \ref{c3} shows. While in $\S$ \ref{characabn}
we do not care about the cost function, in $\S$ \ref{characnorm} the
cost function takes part in the process. To characterize strict
abnormal extremals the cost function is fundamental because these
extremals are abnormal, but not normal. The only way to guarantee
that an extremal is not normal is to use the cost function.

As a consequence of the final constraint submanifolds
obtained from the algorithm for abnormality and normality,
the strict abnormality can be studied. The adjective strict
denotes that the extremal only admits one kind of lifts to
the cotangent bundle. To find strict abnormal extremals we
have to project the final constraint submanifolds to $M$.
In the intersection there are the extremals admitting two
different kinds of lifts: with $p_0=0$ and with $p_0=-1$.
This explanation makes evident the presence of the cost
function to study strict abnormality because the final
constraint submanifold for normality is used.

To sum up, all the biextremals in $N_f^{[0]}$ and
$N_f^{[-1]}$ are projected through $\rho=\pi_M \circ \pi_1
\colon T^*M\times U \rightarrow M$ due to Remark
\ref{nocontrol}.

\begin{prop} Let $P=\rho(N^{[0]}_f) \, \cap \, \rho (N^{[-1]}_f)$.
\begin{itemize}
\item[(i)] If $P=
\emptyset$ and $\rho(N^{[0]}_f)\neq \emptyset$, then all
the abnormal extremals are strict.
\item[(ii)] If $P=
\emptyset$ and $\rho(N^{[-1]}_f)\neq \emptyset$, then all the normal
extremals are strict normal.
\item[(iii)] If $P\neq \emptyset$ and $\rho(N^{[0]}_f)=P$, then there are no strict abnormal
extremals.
\item[(iv)] If $P\neq \emptyset$ and $\rho(N^{[0]}_f) \neq P$, then there are
locally abnormal extremals.
\item[(v)] If $P\neq \emptyset$ and
$\rho(N^{[0]}_f)=\rho(N^{[-1]}_f)=P$, then all the abnormal
extremals are also normal and viceversa.
\end{itemize}
\end{prop}
In item $(iv)$, it is said that there are strict abnormal
extremals, but only locally since the extremal could have
pieces in $P$. So at some points the extremal can be
locally normal.

\section{Free-time optimal control problem}\label{SFOCP}

Once the theory has been introduced let us deal with the
particular case of the free time OCP. In this case the
interval of definition of the extremals is another unknown
of the problem.
\begin{state} \textbf{(Free-time Optimal Control Problem,
$FOCP$)}
\\
Given $M$, $U$, $X$, $F$, $x_a$, $x_b\in M$ (as in $\S$ \ref{SOCP}).
Find $(\gamma,u)$ and $I=[a,b]\subset \mathbb{R}$ such that
\begin{itemize}
\item[(1)] $\gamma(a)=x_a$, $\gamma(b)=x_b$,
\item[(2)] $ \dot{\gamma}(t)=X(\gamma(t), u(t))$, $t\in I$, and
\item[(3)] ${\cal S}[\gamma,u]$ is minimum over all curves on $M\times U$
satisfying  $(1)$ and $(2)$.
\end{itemize}
\end{state}
Pontryagin's Maximum Principle is the same as Theorem \ref{PMPpre},
but replacing $(2.a)$ by
\[(2.a') \; H(\lambda(t),u(t)) \, \sl{ is \,  zero \, everywhere} \, t\in
I.\] Thus the presymplectic equation (\ref{preqP}) must be
restricted to the submanifold defined by the condition
\[H=H_X +p_0F=0.\]
Hence, it must also be stabilized in the algorithm. Due to
the properties of hamiltonian systems
\cite{AbrahamMarsden}, the condition $H=0$ is trivially
stabilized along integral curves of the hamiltonian vector
field. Thus its tangency condition does not add any new
constraint to the submanifolds of the algorithm. The same
happens with $H={\rm constant}$, but this is not a suitable
constraint for a submanifold, that is why it was not
included in the primary constraint submanifold for
fixed-time OCP $\S$ \ref{Scharacterize}. In contrast to
Remark \ref{setextrem}, the final constraint submanifolds
we find here recover exactly the whole set of extremals
since all the necessary conditions of Theorem \ref{PMPpre}
are taking into account. The trivial stabilization of $H=0$
makes possible to run the algorithm putting aside, then the
same final constraint submanifolds as in $\S$
\ref{characabn}, \ref{characnorm} are obtained. Those
submanifolds are renamed, respectively, as $N_{ff}^{[0]}$
and $N_{ff}^{[-1]}$ since the actual final constraint
submanifolds are obtained by considering the vanishing of
the Hamiltonian:
\begin{eqnarray*}
N_f^{[0]}&=&N_{ff}^{[0]}\cap \{(\lambda,u)\in T^*M\times U \, | \, H_X=0\}, \\
N_f^{[-1]}&=&N_{ff}^{[-1]}\cap \{(\lambda,u)\in T^*M\times
U \, | \,H_X-F=0\}.
\end{eqnarray*}
Due to condition $(2.b)$ in Theorem \ref{PMPpre}, the zero fiber
must be deleted from $N_f^{[0]}$.
\begin{prop}\label{time} Given a free-time optimal control
problem:
\begin{enumerate}
\item If $N_f^{[0]}$ has only zero covectors, there are no abnormal
extremals.
\item If $N_f^{[0]}$ has nonzero covectors and $N_{ff}^{[0]}\subset
\{(\lambda,u)\in T^*M\times U \, | \, H_X=0\}$, then every abnormal
extremal is strict and there are no normal extremals as long as $F$
does not vanish.
\end{enumerate}
\end{prop}

\section{Example}\label{SExamples}
There are some classical optimal control problems where the
classification of extremals has been described with
different tools and approaches: geodesics in Riemannian
geometry \cite{LS96}, shortest-paths in subRiemannian
geometry \cite{AgrachevGauthier,LS96} and OCPs with
control-affine systems
\cite{2004Agrachev,AZ,Trelat2000,Trelat2006}. All of them
can be studied in a unified way by direct application of
the method we have proposed in this paper. Here we are
going to use the algorithm for a particular example.

\subsection{Control-affine mechanical system}
Following the described method we find a strict abnormal
extremal for a control-affine system on $TQ$, that, in
fact, models an affine connection control mechanical
system. See more details about these systems in
\cite{2005BulloLewisBook}.

Let $M=TQ=T\mathbb{R}^3$ (i.e. $Q=\mathbb{R}^3$),  $U$ be
an open set in $\mathbb{R}^2$ containing the zero. In local
natural coordinates $(x,y,z,v_x,v_y,v_z)$ for $TQ$, the
drift vector field of the control-system is
\[\ds{Z=v_x\parx+v_y\pary+v_z\parz} \ , \]
and the input vector fields are $\ds{Y_1=\parvx}$ and
$\ds{Y_2=(1-x)\parvy+x^2\parvz}$. So the control system is
given by $Z+u_1Y_1+u_2Y_2$. The endpoint conditions on $TQ$
are $v_a=(2,0,0,0,v_y^0,4(1-v_y^0))$,
$v_b=(2,1,0,0,2(1-v_y^0),4v_y^0-4)$ with $v_y^0\neq 1$. The
cost function is $F=\frac{u_1^2+u_2^2}{2}$. Hence
Pontryagin's Hamiltonian is
\[H(\lambda,u_1,u_2)=p_xv_x+p_y v_y +p_z v_z+u_1q_x+u_2(1-x)q_y+u_2x^2q_z+p_0\frac{u_1^2+u_2^2}{2}\]
with Hamilton's equations for abnormality and normality
\[\begin{array}{llll}
\dot{x}=v_x & \dot{v_x}=u_1
 & \dot{p}_x=q_y u_2-2q_zu_2 x & \dot{q}_x=-p_x\\
\dot{y}=v_y & \dot{v_y}=u_2 (1-x)
  & \dot{p}_y=0 &\dot{q}_y=-p_y\\
\dot{z}=v_z & \dot{v_z}=u_2 x^2
  & \dot{p}_z=0 &
\dot{q}_z=-p_z
\end{array}\]
and Hamiltonian vector field $X_H=\sum_{i\in\{x,y,z\}}
(A^i\partial/ \partial i +B^i\partial/ \partial
v_i+C_i\partial/ \partial p_i+D_i\partial/ \partial
q_i)+E_1\partial/ \partial u_1+E_2\partial/ \partial u_2$,
where $A^i, B^i, C_i, D_i$ are determined by Hamilton's
equations.

The constraint algorithm for abnormality gives us
\[\begin{array}{lclcl}
N_0^{[0]}&=&\{(\lambda,u) \in T^*TQ\times U & | &
\partial H^{[0]} /\partial u_k(\lambda,u)=H_{Y_k}(\lambda)=0  \; , \; k=1,2\} \\
&=&\{(\lambda,u) \in T^*TQ\times U & | & q_x=0, \;
\mathbf{q_y(1-x)+q_zx^2=0 }\}\\N_1^{[0]}&=&\{(\lambda,u)
\in N_0^{[0]}  &|&
H_{[Z,Y_k]}(\lambda)=0 \quad {\rm for} \;\; k=1,2 \} \\
&=&\{(\lambda,u) \in N_0^{[0]} & |& p_x=0, \;
(-1+x)p_y-x^2p_z-v_xq_y+2xv_xq_z=0\}
\\ N_2^{[0]}&=&\{(\lambda,u) \in N_1^{[0]} &|&
(H_{[Z,[Z,Y_k]]+u_l[Y_l,[Z,Y_k]]})(\lambda)=0 \quad {\rm
for} \;\; k=1,2\}\\
&=&\{ (\lambda,u) \in N_1^{[0]}  & | &
\mathbf{(-q_y+2xq_z)u_2=0}, \\ &&&&
(-q_y+2xq_z)u_1=-(2p_yv_x-4xv_xp_z+2v_x^2q_z) \} \
.\end{array}\] In order to satisfy the endpoint conditions,
not to have the zero covector and to have a strict abnormal
extremal, the subset defined by $x \, (x-1) \, q_z \,
u_2=0$, coming from the above bold equations, must be
deleted from the constraint submanifold. Then
\[\begin{array}{lcl}
N_2^{[0]}&=&\{(\lambda,u)\in  T^*TQ\times U-\{x(x-1)q_zu_2=0\}\, | \, q_x=0, \; -q_y+4q_z=0, \; p_x=0, \\
&&
p_y-4p_z=0, \; x=2, \; v_x=0\} \\
N_3^{[0]}&=&\{(\lambda,u) \in N_2^{[0]} \, | \, v_x=0, \;
u_1=0\}\\
N_4^{[0]}&=&\{(\lambda,u) \in N_2^{[0]} \, | \, u_1=0, \;
E_1=0\}=N_5^{[0]}=N_f^{[0]}.
\end{array}\]
By restriction to the final constraint submanifold and
integrating Hamilton's equations on $[0,1]$ we have
$\lambda(t)=(0,4p_z^0,p_z^0,0,-4p_z^0t+4q_z^0,-p_z^0t+q_z^0)$
and
\[\gamma(t)=(2,-u_2\frac{t^2}{2}+v_y^0t,2u_2t^2+4(1-v_y^0)t,0,-u_2t+v_y^0,4u_2t+4(1-v_y^0)) \, \]
with $u_2=2(v_y^0-1)$.

The constraint algorithm for normality gives us
\[\begin{array}{lclcl}
N_0^{[-1]}&=&\{(\lambda,u) \in T^*TQ\times U &| &
\partial H^{[-1]} /\partial u_k(\lambda,u)=0  \; , \; k=1,2\} \\
&=&\{(\lambda,u) \in T^*TQ\times U &| & q_x-u_1=0 \ , \\
&&&& q_y (1-x)+q_zx^2-u_2=0\}=N_1^{[-1]}=N_f^{[-1]}.
\end{array}
\] If we substitute the curve $\gamma$ in Hamilton's
equations, we have $u_1=0$ and $u_2=2(1-v_y^0)$, then for
the primary constraint submanifold $q_x=0$ and
$u_2=q_y-4q_z$. Due to Hamilton's equations $p_x=0$ and
$0=\dot{p}_x=u_2^2$. This last equality is only possible if
$v_y^0=1$, but that was not the hypothesis. Thus there does
not exist a covector with $p_0=-1$ along $\gamma$, that is,
$\gamma$ is a strict abnormal extremal whenever $v_y^0\neq
1$.

\section{Conclusion and outlook}

In this paper we have given a geometric method to study
different kinds of extremals in a wide range of optimal
control problems with an open control set. This can also be
applied in the case of closed control set following ideas
in \cite{Lopez}. This method is based on the suitable
reinterpretation of the so-called {\sl presymplectic
algorithm} in other fields. The dependence on the cost
function makes difficult to give general characterizations
of normal and strict abnormal extremals since each problem
must be studied by itself. However, the abnormal extremals
only depend on the geometry of the control system, so some
general results can be deduced. See \cite{1997BonnardKupka}
for an approach to a related problem for single-input
control-affine systems.

One line of future research is to apply this general
algorithm in the study of optimal control problems with
affine connection control systems, which model the motion
of different types of mechanical systems such as rigid
bodies, nonholonomic systems and robotic arms
\cite{2005BulloLewisBook}. Eventually, we will focus on
particular problems for mechanical systems, as for instance
time-optimal problems and control-quadratic cost function.

Apart from having sufficient conditions to determine where
the extremals are, it may be interesting to prove the
density and the optimality of them, similar to the work
done in \cite{LS96} for control-linear systems with two
input vector fields. Here the cost function takes part in,
even in determining the optimality of abnormal extremals.
There are some results that do not contribute to the
optimism in relation to the existence of strict abnormal
minimizers, at least in a generic sense. For instance, in
\cite{ChitourJeanTrelatDiffGeom2006,ChitourJeanTrelat2006}
it is proved the existence of an open and dense subset in
the set of control systems where every nontrivial strict
abnormal extremal is not a minimizer for control-quadratic
cost functions and control-affine systems.

Another meaningful research issue is to establish
connections between the controllability of the system
\cite{2005BulloLewisBook,VanderSchaft,SJ72} and the final
constraint submanifold obtained for abnormality, since both
notions are independent of the functional to be minimized.
In fact, we are already working on some properties of
controllability, similar to results in \cite{BS93}, that
can also be justified using the algorithm here described.

\section*{Acknowledgements}
We acknowledge the financial support of \emph{Ministerio de
Educaci\'on y Ciencia}, Project MTM2005-04947 and the
Network Project MTM2006-27467-E/. MBL also acknowledges the
financial support of the FPU grant AP20040096. We thank
Professor Andrew D. Lewis and Professor David Mart\'in de
Diego for their useful ideas to develop this work.

\end{document}